\documentclass[12pt]{amsart}

\usepackage{latexsym}
\usepackage{amsmath}
\usepackage{amssymb}
\usepackage{amscd}
\usepackage{multicol}
\usepackage[cmtip,matrix,arrow]{xy}

\newtheorem{theorem}{\bf Theorem}[section]

\newtheorem{proposition}[theorem]{\bf Proposition}
\newtheorem{lemma}[theorem]{\bf Lemma}

\newtheorem{definition-theorem}[theorem]{\bf Theorem-Definition}

\def\Q{\mathbb{Q}}

\def\R{\mathbb{R}}
\def\bC{\mathbb{C}}

\def\g{\mathfrak{g}}
 \def\k{\mathfrak{k}}

\def\g{\frak{g}}

\def\p{\frak{p}}
\def\a{\frak{a}}

\def\k{\frak{k}}

\def\H{\mathbb{H}}
\title[Connectivity  for isoparametric submanifolds]{Connectivity and Kirwan surjectivity
 for isoparametric submanifolds}

\date{\today}
\author[A.-L. Mare]{Augustin-Liviu  Mare}
\address{
Department of Mathematics and Statistics\\ University of
Regina \\ College West 307.14 \\ Regina SK, Canada S4S
0A2}
\email{mareal@math.uregina.ca}

\begin{document}
\begin{abstract}
Atiyah's formulation of what is nowadays called the  convexity theorem of Atiyah-Guillemin-Sternberg has
two parts: (a) the image of the moment map arising from a Hamiltonian action of a torus on a symplectic manifold
is a convex polytope, and (b) all preimages of the moment map are connected. Part (a) was generalized
by Terng to
the wider context of isoparametric submanifolds in euclidean space. In this paper we  prove a
generalization of
part (b) for a certain class of isoparametric submanifolds (more precisely, for those with all multiplicities
strictly greater than 1). For generalized real flag manifolds,
which are an important class of isoparametric submanifolds, we  give a surjectivity criterium
for a certain Kirwan map (involving equivariant cohomology with coefficients in $\Q$) which arises naturally
in this context. Examples are also discussed.  

\end{abstract}
\maketitle

\section{Introduction}

A crucial result in symplectic geometry is the convexity theorem of Atiyah-Guillemin-Sternberg.
It states that if $(M, \omega)$ is a symplectic manifold acted on by a torus $T$ in a Hamiltonian
way, then the image of the moment map $\mu:M \to {\rm Lie}(T)^*$ is a convex polytope. 
A closely related result says [A] that all pre-images of $\mu$ are connected (or empty). 
The convexity theorem has been generalized by Terng [T1] to the case when $M$ is an
isoparametric submanifold in a euclidean space (see section 1 of our paper for the definition;
we  note that we always assume  that $M$ is compact). 
More precisely, she proved the following theorem.
\begin{theorem}\label{terng} {\rm (Terng [T1])} Let  $M\subset \R^{n+k}$ be  an isoparametric submanifold,
$q$ a point of $M$,  $\nu_q(M)$ the normal space at $q$ to $M$, and $P:\R^{n+k}\to \nu_q(M)$
the orthogonal projection map.  Then
 the image of the  map $$\mu=P|_{M}:M \to \nu_q(M)$$  is a convex polytope. 
\end{theorem}
For example, $M$ can be a generalized real flag manifold, i.e. an orbit of
the isotropy representation of a compact symmetric space (for more details, see section
\ref{kirwan} of our paper). These can be realized
as real loci (i.e. fixed point sets of antisymplectic involutions) 
of coadjoint orbits. Convexity results for such spaces had been obtained by Kostant [K] (see
also Duistermaat [D]). But Terng's theorem does not fit entirely into the framework of the paper
by Kostant.
The reason is that there are several examples of isoparametric submanifolds
which are not flag manifolds, i.e. are not homogeneous: such examples were found by Ferus, Karcher and
M\"unzner in [FKM]. 

Now turning to the connectedness result of Atiyah mentioned above, the following question
can be raised. In the context of Theorem \ref{terng}, is any nonempty preimage of $\mu$ connected?
One can easily see that the answer is in general negative. More specifically, the real flag
manifold corresponding to the symmetric space $SU(2)/SO(2)$ is a circle and $\mu$ is just the height 
function on that  circle; hence almost all preimages  consist of two points.  
Our main result is (see section 2 for the definition of multiplicities):
\begin{theorem}\label{main} Let $M\subset \R^{n+k}$ be an isoparametric submanifold
with all multiplicities   greater or equal to 2, and let $a\in \nu_q(M)$ be in the image
of $\mu$.  The following statements are true.

(i) The function $f:M \to \R$, $f(x) =\| \mu(x)-a\|^2$, is minimally degenerate
in the sense defined by Kirwan {\rm [K]} (see section 3 of our paper).

(ii) The level set $\mu^{-1}(a)$ is connected or empty.
\end{theorem}
 
 We note that in the case of the real flag manifold arising from $SU(2)/SO(2)$ mentioned above, 
the (only) multiplicity is equal  to 1. So the hypothesis on multiplicities in the previous theorem is 
essential.

\noindent {\bf Remarks.}
(a) In fact, Terng [T1] proved that if
$M_{\xi}$ is a  manifold parallel to $M$ (see section 2), then the image of
$P|_{M_{\xi}} : M_{\xi} \to \nu_q(M)$ is a convex polytope. With the methods of our paper one can prove that   
any preimage of the latter map is connected as well.

(b) Here is a list of examples of isoparametric submanifolds with all multiplicities
at least equal to 2, together with the corresponding
 parallel manifolds (see the previous remark):
\begin{itemize}
\item adjoint orbits of compact Lie groups
\item isotropy orbits of the symmetric spaces $SU(2n)/Sp(n)$, $E_6/F_4$, $SU(m+n)/S(U(m)\times U(n))$,
where $m>n$,
$Sp(m+n)/Sp(m)\times Sp(n)$ (for more details, see Helgason [H])
\item infinitely many of the  Ferus-Karcher-M\"unzner [FKM] examples. 
\end{itemize} 

(c) The idea of the proof of Theorem \ref{main} goes back to Kirwan [K]. We consider the  Morse stratification
of $M$ induced by $f$. The strata  are smooth submanifolds, and all but the minimum one have
codimension at least 2,   due to the hypothesis on multiplicities.
Hence the minimum stratum is connected.
Finally, we use the fact, also proved by Kirwan [K], that the latter stratum has the same 
$\check{\rm C}$ech cohomology module
$H^0(\cdot, \Q)$
as the minimum set $f^{-1}(0)=\mu^{-1}(a)$.

(d) Terng [T2] was able to extend Theorem \ref{terng} to the case when $M$ is an {\it infinite dimensional} 
isoparametric submanifold in a Hilbert space, or any of its parallel manifolds, like in Remark (a).
We conjecture that also our Theorem \ref{main} can be extended to the infinite dimensional situation. 
The connectivity theorem proved recently  by Harada, Holm, Jeffrey, and the author in [HHJM]
shows that the conjecture is true for the loop group $\Omega(G)$, which is a parallel manifold of a 
certain infinite dimensional isoparametric submanifold with all multiplicities equal to 2.

\noindent {\bf Acknowledgement.} I would like to thank Jost Eschenburg for discussions and important hints
concerning the topics of the paper, especially those in the last section. 

\section{Isoparametric submanifolds in euclidean space}

We  present some basic facts concerning the theory of isoparametric submanifolds.
For more details, the  reader can consult  Palais and Terng's book [PT]
(especially chapter 6) and the references therein.

Let $M\subset \R^{n+k}$ be an $n$-dimensional
embedded submanifold, which is  closed, complete with respect to the induced metric,
and full
(i.e. not contained in any affine subspace). We say that $M$ is
  {\it isoparametric}  if  any normal vector at a point
of $M$
 can be extended
 to a parallel normal vector field $\xi$ on $M$ with the property that the eigenvalues of the shape operators
$A_{\xi(p)}$ (i.e. the principal curvatures) are independent on $p\in M$, as values and multiplicities.
It follows that for  $p\in M$, the set  $\{A_{\xi(p)}: \xi(p) \in \nu (M)_p\}$ is a commutative family of
selfadjoint endomorphisms of $T_p(M)$, and so it determines a decomposition of $T_p(M)$  as a
direct sum of common eigenspaces
$E_1(p), E_2(p), ...,E_r(p)$. There exist normal vectors  $\eta_1(p), \eta_2(p),..., \eta_r(p)$
such that
$$A_{\xi(p)}|_{E_i(p)}=\langle \xi(p), \eta_i(p) \rangle {\rm id}_{E_i(p)},$$
for all $\xi(p) \in \nu_p (M)$, $1\le i\le r$. By parallel extension in the normal bundle we obtain
the vector fields $\eta_1, \ldots, \eta_r$, which are the {\it principal curvature vectors}.
The eigenspaces from above give rise to the distributions $E_1, \ldots, E_r$ on $M$, which are called
the  {\it curvature distributions}.
The numbers $$m_i= {\rm rank} E_i,$$ $1\leq i \leq r$, are  the {\em
multiplicities } of $M$.

To any parallel normal vector field $\xi$ on $M$ we assign the end-point map
$\pi_{\xi}:M \to \R^{n+k}$, $$\pi_{\xi}(p) = p+\xi(p),$$ $p\in M$. Its image is the
``parallel" manifold $M_{\xi}$, which is also an embedded submanifold of $\R^{n+k}$.
Now the differential map of $\pi_{\xi}$ is
$$d(\pi_{\xi})_p= I - A_{\xi(p)}.$$
So the focal points of $M$ in the affine normal space $p+\nu_p(M)$ are those which are in
one of the hyperplanes
$$\ell_i(p):= \{p+\xi(p) : \langle \eta_i(p), \xi(p) \rangle =1\},$$
for some $i\in \{1,2,\ldots, r\}$. It turns out that
$\ell_1(p),\ldots, \ell_r(p)$ have a unique intersection point, call
it $c_0$, independent on the point $p$. Moreover, $M$ is contained
in a sphere centered at $c_0$ (here we use the assumption that $M$ is compact). 
We do not lose any generality if we
assume that $M$ is contained in the unit sphere $S^{n+k-1}$, hence
$c_0$ is just the origin $0$. One shows that the group of linear
transformations of $\nu_p(M)$ generated by the reflections about
$\ell_1(p),\ldots, \ell_r(p)$ is a Coxeter group, whose isomorphism
type is  independent on $p$. We denote it by $W$ and call it the
{\it Weyl group} of $M$.

The map $\pi_{\xi} : M \to M_{\xi}$ is a submersion.
If $p$ is a point in  $M$ and $b:=\pi_{\xi}(p)=p+\xi(p)$, we denote by
$S_{p, b}$ the connected component of $\pi_{\xi}^{-1}(b)$ which contains $p$.
This manifold is called  the {\it slice} through $p$ corresponding to $\xi$. Consider
$$I:=\{i\in \{1,2,\ldots, r\} : p+\xi(p) \in \ell_i(p)\},$$
and also the subspace $(\bigcap_{i\in I} \ell_i(p))^{\perp}$ of $\nu_p(M)$, where the
subscript $\perp$ indicates the orthogonal complement. An important result
is the so-called slice theorem, which is stated as follows (for more details, see
[PT, Theorem 6.5.9]).
\begin{theorem} The slice $S_{p, b}$ is a (full) isoparametric submanifold of $(\bigcap_{i\in I} \ell_i(p))^{\perp}
+\sum_{i\in I}E_i(p)$.
\end{theorem}

\section{Minimally degenerate functions according to Kirwan}

In this section we follow chapter 10 of Kirwan's book [K].

\vspace{0.5cm}

\noindent {\bf Definition.}  {\it A smooth function $f: X \to \R$ on a closed manifold $X$ is
called {\rm minimally degenerate} if the following conditions hold.
\begin{itemize}
\item[(a)] The set of critical points of $f$ is a finite union of disjoint  closed subsets
$C_1, C_2, \ldots, C_N$, on each of which $f$ is  constant.   These subsets  are called the
{\rm critical sets} of $f$.
\item[(b)] For every $j=1, \ldots, N$, there exists a submanifold $Y_j$ of $X$, which contains
$C_j$, with the property that the normal bundle\footnote{More specifically,
$(TX / TY_j)|_{Y_j}$} of $Y_j$ in $X$ is orientable, and such that
\begin{itemize}
\item[(i)] $C_j$ is the subset of $Y_j$ on which $f$ takes its minimum value
\item[(ii)] for any $x\in C_j$, the space $T_xY_j$ is maximal among the subspaces of
$T_xX$ on which the Hessian ${\rm Hess}_x(f)$ is positive semi-definite.
\end{itemize}
A manifold $Y_j$ with the properties (i) and (ii) is called a {\rm minimizing manifold}
around $C_j$. The codimension of $Y_j$ is called the {\rm index} of $f$ along $C_j$.
\end{itemize}
}

Even though minimally degeneracy is a condition weaker than nondegeneracy in the sense of Bott,
it is still sufficient to induce a Morse stratification of $M$, as the following theorem shows.
\begin{theorem}\label{stratification} {\rm (Kirwan [Ki])} Let $f: X \to \R$ be a minimally degenerate function
like above and let $g$ be a Riemannian metric on $M$. For any $j\in \{1,2,\ldots, N\}$ we denote by
$\Sigma_j=\Sigma_j(g)$ the set of all points in $M$ with the property that the $\omega$-limit of the integral line
through $x$
of the vector field $-\nabla (f)$ is contained in $C_j$. Then we have as follows.

(a) There exists a metric $g$ with the
property that $\Sigma_j$ is a smooth submanifold of $M$ of codimension
equal to the index along $C_j$. We also have
$$M=\bigcup_{1\le j\le N} \Sigma_j, \quad \Sigma_i\cap \Sigma_j =\phi, \ {\it for \ } i\neq j.$$
The intersection of  $Y_j$ with a sufficiently small neighbourhood
of $C_j$ is contained in $\Sigma_j$.  

(b) For any $j\in \{1,2,\ldots, N\}$, the inclusion map $C_j\hookrightarrow \Sigma_j$ induces an isomorphism
in  $\check{\it C}$ech cohomology.
\end{theorem}

\section{Proof of Theorem \ref{main}}\label{four}

Let $M^n\subset \R^{n+k}$ be an isoparametric submanifold.
We start  with the following lemma.
\begin{lemma}\label{lem} Let $q$ be a point of $M$, $b \in \nu_q(M)$, and let $S_{q,b}$ be  the corresponding
slice. Consider $$I:=\{i\in \{1,2,\ldots, r\} :
b\in \ell_i(q)\}.$$
If $x$ is an arbitrary point in $S_{q,b}$, then we have:

(a) $S_{q,b}$ is a full isoparametric submanifold in $(\bigcap_{i\in I} \ell_i(q))^{\perp}
+\sum_{i\in I}E_i(q)$, whose normal space at  $q$ is $(\bigcap_{i\in I} \ell_i(q))^{\perp}$

(b) $T_x S_{q,b}= \sum_{i\in I}E_i(x)$

(c) $\bigcap_{i\in I} \ell_i(x)+ \sum_{i\in I} E_i(x)  = \bigcap_{i\in I} \ell_i(q)+\sum_{i\in I} E_i(q) $

(d)  $\bigcap_{i\in I} \ell_i(q)\subset \nu_x(M) \cap \nu_{q}(M)$

(e) $\sum_{j\notin I} E_j(x)$ is perpendicular to $\nu_{q}(M)$.

\end{lemma}

\begin{proof} Points (a)-(d) have been proved for instance in [PT, Chapter 6].
We will prove (e). From (c) we deduce
that $\sum_{j\notin I} E_j(x)$ is perpendicular to $\bigcap_{i\in I} \ell_i(q)$, and from (d), the same
space is perpendicular to $\bigcap_{i\in I} \ell_i(q)$.

\end{proof}

We study the function
 $f: M \to \R$, $f(x) = \|\mu(x)-a \|^2$.
 First we determine its critical points. We note that this has been done already by Terng in
[T2, section 3]. Let us consider an orthonormal
basis $e_1, \ldots, e_k$ in $\nu_q(M)$. We have $$\mu(x)-a= \sum_{i=1}^k \langle x-a, e_i\rangle e_i,$$
hence
$$f(x) =  \sum_{i=1}^k \langle x-a , e_i\rangle ^2,$$
which implies \begin{equation}\label{gradient} df_x(v)= 2\sum_{i=1}^k \langle x-a, e_i\rangle \langle v, e_i\rangle = 2\langle v,
 \sum_{i=1}^k \langle x-a, e_i\rangle e_i \rangle = 2\langle
v, \mu(x)-a \rangle,\end{equation}
for any $v\in T_x(M)$.  So $x$ is a critical point of $f$ if and only if the vector $b:= \mu(x)-a\in \nu_q(M)$
is perpendicular to
$T_x(M)$, in other words, when $x$ is a critical point of the height function
$h_b : M\to \R$, $h_b(x)=\langle b,x \rangle$. We can express this in a more concise way as
\begin{equation}\label{critical}Crit(f)= \mu^{-1}(a) \cup \bigcup_{b\in \nu_q(M)} \mu^{-1}(a+b)\cap Crit(h_b).\end{equation}
We  prove that the intersection in the right hand side is nonempty  only for finitely many $b\in \nu_q(M)$.
We  use the fact  that $$Crit(h_b) =\bigcup_{w\in W} S_{wq,b}$$ (see e.g. [T2, subsection 2.9]).
\begin{lemma} There exists finitely many $b\in \nu_q(M)\setminus \{0\}$ with the property that
$\mu^{-1}(a+b)\cap S_{wq,b}$ is nonempty.
\end{lemma}
\begin{proof} For simplicity we assume that $w=1$.
Pick $I$ an arbitrary subset of $\{1,2,\ldots, r\}$.
In $\nu_q(M)$ we consider the subspace
$$\ell_I:=\bigcap_{i\in I} \ell_i(q).$$
We also consider the convex hull  ${\rm cvx} [(W_I). q]$, where $W_I$ denotes the
stabilizer if $\ell_I$.  The affine span of this convex body  is  perpendicular to $\ell_I$ and
its dimension is just the codimension of $\ell_I$. To justify this, we note that this affine span
  is  the
affine normal space to a certain slice, which has codimension equal to $\dim (\ell_I)^{\perp}$,
and which is contained in an affine space perpendicular to $\ell_I$ (see Lemma \ref{lem} (a)).
    Consequently the intersection
$(a+\ell_I)\cap {\rm cvx} [(W_I). q]$ has at most one point. By letting
$I$ vary among all subsets of $\{1,2,\ldots, r\}$, we obtain a finite set of points, call it $F$.
Now, if $\mu^{-1}(a+b)\cap S_{wq,b}\neq \phi$, then
$a+b$ must belong to $\mu(S_{wq,b})$, which is ${\rm cvx}[(W_b).q]$. Consequently, we have
$$a+b \in (a+\ell_I)\cap {\rm cvx} [(W_I). q],$$
where $I=\{i\in \{1,2,\ldots, r\} : b\in \ell_i(q)\}$. This implies $a+b\in F$.
\end{proof}

We have proved that the set
 $$B:=\{b\in \nu_q(M)\setminus \{0\} : \mu^{-1}(a+b)\cap Crit(h_b)\neq \phi\}$$
is finite.  The description (\ref{critical}) can be refined by taking into account that for
any $b\in \nu_q(M)$ we have $h_b(x) = \langle  x,b \rangle = \langle \mu(x), b\rangle$, $x\in M$.
Consequently, if $x \in \mu^{-1}(a+b)$, then $h_b(x) = \langle a+b, b\rangle$. So
\begin{equation}\label{refined}Crit(f) = \mu^{-1}(a)\cup \bigcup C_{b,w}\end{equation}
where $C_{b,w}:=\mu^{-1}(a+b)\cap S_{wq,b}$ and the union runs
over all $b\in B$ and  $w\in  W$  with the property  that $h_b(S_{wq,b})=\langle a+b, b\rangle$.

Fix $b$ and $w$ like above. Let $Y_{b,w}$ denote the stable manifold of the function $h_b$ corresponding
to the critical set $S_{wq,b}$. More specifically, this consists of all points $x\in M$ with the property
that the limit at $\infty$ of the integral line through $x$ of the vector field $-\nabla(h_b)$ is
in $S_{wq,b}$. We will prove the following result.
\begin{proposition}\label{minimal} (a) A minimizing manifold for $f$ around $\mu^{-1}(a)$ is $M$ itself.

(b) For $b\in B$ and $w\in W$, the space $Y_{b,w}$ is a minimizing manifold for $f$ around $C_{b,w}$.

(c) The only critical set of index 0 is $\mu^{-1}(a)$.
\end{proposition}

\begin{proof}
(b) First we show that the normal bundle to $Y_{b,w}$ is orientable. To this end we note
that $Y_{b,w}$ is just a vector bundle over $S_{wq, b}$, hence it is homeomorphic to the
latter space. But $S_{wq,b}$ is an isoparametric submanifold with all multiplicities at least 2,
hence it is simply connected. Consequently, $Y_{b,w}$ is also simply connected.
This implies the desired conclusion (we recall that any vector bundle over a simply connected
space is orientable).

Next we note that
if $x\in Y_{b,w}$, then $h_b(x) \ge h_b(S_{wq,b}) = \langle a+b, b\rangle.$ So we have
$\langle \mu(x), b\rangle \ge \langle a+b, b \rangle$, which implies $$\langle \mu(x)-a, b\rangle \ge
\langle b,b\rangle.$$
We deduce
\begin{equation}\label{ineq}\|\mu(x)-a\| \cdot \|b\| \ge \langle \mu(x)-a, b\rangle \ge \|b\|^2,\end{equation}
hence
$$f(x) \ge f(C_{b,w}).$$
Moreover, if $x\in Y_{b,w}$ has the property that $f(x) = f(C_{b,w})$, then we must have
\begin{itemize}
\item $h_b(x) = h_b(S_{wq,b})$ (from equation (\ref{ineq})), hence $x\in S_{wq,b}$
\item $\mu(x)-a=\lambda b$, for a number $\lambda$ (from equation (\ref{ineq})); we deduce that $\lambda=1$,
because $\langle \mu(x)-a, b\rangle = \langle b, b\rangle$.
\end{itemize}
Consequently, $x\in C_{b,w}$.
We have proved that the condition (b) (i) from the definition of a minimally degenerate function is satisfied.

It remains to check condition (b) (ii).
Let us consider  a  point $x_0$ in $C_{b,w}$. We  construct a subspace $V\subset T_{x_0}(M)$ with the
following properties.
\begin{itemize}
\item[1.] $V\oplus T_{x_0}(Y_{b,w}) = T_{x_0}(M)$
\item[2.] ${\rm Hess}_{x_0}(f)|_V$ is negative definite.
\end{itemize}
First we determine the Hessians of $h_b$ and $f$ at the point $x_0$.
To this end we consider the functions $H_b$ and $F$  given by
$$H_b(x) = \langle x, b\rangle,\quad F(x) =\| P(x)-a\|^2,\quad x\in \R^{n+k}$$
where $P$ denotes the orthogonal projection $\R^{n+k}\to \nu_q(M)$.
We know that for $v,w\in T_{x_0}(M)$ we have
$${\rm Hess}(f)_{x_0}(v,w) :=\langle \partial_v (\nabla f)(x_0), w\rangle =\langle  \partial_v (\nabla F)(x_0), w\rangle
+\langle A_{(\nabla F)^{\perp}_{x_0}}v,w\rangle,$$
where $\nabla$ stands for gradient and  the superscript $\perp$ indicates the orthogonal projection on $\nu_{x_0}(M)$.
Because $(\nabla F)_{x}= 2(P(x) -a)$ (see equation (\ref{gradient})), and $P(x_0)-a=b$, we deduce that
\begin{equation}\label{hessf}{\rm Hess}(f)_{x_0} = 2(P + A_{b}).\end{equation}
Similarly, the Hessian of $h_b$ is
\begin{equation}\label{hessh}{\rm Hess}(h_b)_{x_0} = A_b.\end{equation}
The tangent space $T_{x_0}(Y_{b,w})$ is the subspace of $T_{x_0}(M)$ where the hessian ${\rm Hess}(h_b)_{x_0}$
is negative semidefinite.
From equation (\ref{hessh}) we deduce that
 for $i\in \{1,2,\ldots, r\}$, the restriction of ${\rm Hess}(h_b)_{x_0}$ to $E_i(x_0)$ is given by scalar
multiplication by  $\langle b, \eta_i(x_0))\rangle$. Consequently we have $$T_{x_0}(Y_{b,w})=\sum E_i(x_0)$$ where the sum
runs over all $i$ with $\langle b, \eta_i(x_0)\rangle \ge 0$. From elementary Morse theoretical
considerations, we know  that
$$T_{x_0}(S_{wq,b})\subset T_{x_0}(Y_{b,w}).$$
Set
\begin{equation}\label{v} V:=\sum E_i(x_0),\end{equation}
where the sum runs over all $i$ with $\langle b, \eta_i(x_0)\rangle <
0$. By Lemma \ref{lem} (e) (with $q$ replaced by $wq$), $P$ maps $V$ to 0. The reason is that
$b\in \ell_i(x_0)$ exactly when
$\langle b, \eta_i(x_0)\rangle =0$ (since $\ell_i(x_0)$ goes through the origin).
By (\ref{hessf}) and (\ref{hessh}), the restrictions
of ${\rm Hess}(f)_{x_0}$ and ${\rm Hess}(h_b)_{x_0}$ to $V$ are the same, 
up to a factor of 2. Since the latter restriction
is strictly negative definite, the former is like that as well.

(c) Suppose that for $b\in B$ and $w\in W$ the critical set $C_{wq,b}$ has index 0. The considerations from
above show that the index of $h_b$ at $S_{wq,b}$ is 0. This implies that the slice $S_{wq,b}$ is the minimum set of
$h_b$ on $M$. Consequently, for any $x\in M$ we must have $h_b(x) \ge h_b(S_{wq,b})$. Like in the proof
of point (b), this implies $f(x)\ge f(C_{b,w})=\|b\|^2$, for all $x\in M$. This is a contradiction
(because $f$ can reach the value 0), which concludes the proof.
\end{proof}

\noindent {\bf Remark.} The arguments used above are approximately the same as  those used by Kirwan
 in the proof of [K, Proposition 4.15] (in the context of Hamiltonian torus actions on symplectic
manifolds).

{\it Proof of Theorem \ref{main}} Only point (ii) has to be proved.
By  Theorem \ref{stratification}, there exists a metric $g$ on $M$ which induces the stratification
$$M= \Sigma_0 \cup \bigcup_{b\in B, w\in W} \Sigma_{b,w},$$
such that $\Sigma_0, \Sigma_{b,w}$ are smooth submanifolds. By Proposition \ref{minimal} and equation (\ref{v}),
the manifolds $\Sigma_{b,w}$ have codimension at least 2 (recall that, by hypothesis, we have
$\dim (E_i(x_0)) = m_i\ge  2$). Consequently, $\Sigma_0$ is connected.
Since the map $\mu^{-1}(a) \hookrightarrow \Sigma_0$ induces the linear
isomorphism $H^0(\mu^{-1}(a), \Q)\simeq H^0(\Sigma_0, \Q)$ (see Theorem \ref{stratification}),
we deduce that $\mu^{-1}(a)$ is connected as well.
 \hfill $\square$

\section{Kirwan surjectivity for real flag manifolds}\label{kirwan}

An important class of isoparametric submanifolds are the real flag manifolds.
More specifically, we start with   a non-compact symmetric space $G/K$, where $G$ is a
non-compact connected semisimple Lie group and $K\subset G$ a
maximal compact subgroup. Then $K$ is the fixed point set of an
automorphism $\tau$ of $G$ (see for instance [H,  chapter
VI]).  The differential map $d(\tau)_e$ is an automorphism of  $\g={\rm Lie}(G)$
and  induces the Cartan
decomposition
$$\g =\k \oplus \p,$$ where $\k={\rm Lie}(K)$ and $\p$ are the
$(+1)$-, respectively $(-1)$-eigenspaces of $(d\tau)_e$.

Now let us consider $\a \subset \p$ a maximal abelian subspace. The number
$k:=\dim(a)$ is the rank of the symmetric space $G/K$.  The
roots of the symmetric space  are linear functions $\alpha
:\a\to \R$ with the property that the space
$$\g_{\alpha}:=\{z\in \g : [x,z]=\alpha(x)z \ {\rm for \ all }
\ x\in\a\}$$ is non-zero. The set $\Pi$ of all roots is a root
system in $(\a^*, \langle \ , \ \rangle)$. Let $\Pi^+\subset \Pi$ be the
set of positive roots with respect to a simple root system. For any $\alpha\in\Pi^+$ we
have
$$\g_{\alpha}+\g_{-\alpha}=\k_{\alpha}+\p_{\alpha},$$
where $\k_{\alpha}=(\g_{\alpha}+\g_{-\alpha})\cap\k$ and
$\p_{\alpha}=(\g_{\alpha}+\g_{-\alpha})\cap\p$. We have the direct
decompositions
$$\p=\a+\sum_{\alpha \in \Pi^+}\p_{\alpha},\quad \k=\k_0+\sum_{\alpha \in \Pi^+}\k_{\alpha},$$
where $\k_0$ denotes the commutator of $\a$ in $\k$. 

Since
$[\k,\p] \subset \p$, the space $\p$ is $Ad_G(K):=Ad(K)$-invariant.
The orbits of the action of $Ad(K)$ on $\p$ are called {\it generalized real
flag manifolds}. The restriction of the Killing
form of $\g$ to $\p$ is an $Ad(K)$-invariant inner product on $\p$,
which we denote  by $\langle \ , \ \rangle$.

\begin{proposition}{\rm (see e.g. [PT, Example 6.5.6])} Let $M=Ad(K).q$ be the orbit of $q\in \a$.

(i) If $q$ is regular (i.e. not contained in any of the hyperplanes $\ker (\alpha)$,
$\alpha \in \Pi$), then $M$ is an isoparametric submanifold of $(\p, \langle \ , \ \rangle)$.
The curvature distributions at $q$ are
$$E_{\alpha}(q) = [q, \k_{\alpha}+\k_{2\alpha}],$$
where $\alpha$ is a positive indivisible root. Hence the multiplicities of $M$ are
$$m_{\alpha}= \dim(\k_{\alpha})+\dim(\k_{2\alpha}).$$
The normal space to $M$ at the point $q$ is $$\nu_q(M)=\a.$$
   
   (ii) If $q$ is not regular, then $M$ is a manifold parallel to an isoparametric submanifold.
   \end{proposition}
   
In [M] we have investigated the action of 
$$K_0:= Z_K(\a) =\{k\in K : Ad(k)(x)=x, \ {\rm for\ all \ } x\in \a\}$$
on $M$.
It turns out that  if all $m_{\alpha}$ are strictly greater than 1, then $K_0$ is connected
and the action of $K_0$ on $M$ is equivariantly formal. In this paper  we will address the following 
question.

\noindent {\bf Problem.} {\it If $\mu: M\to \a$ is the restriction to $M$ of the orthogonal projection map 
$P: \p \to \a$ and $a$ is an arbitrary point in $\a$, is it true that the Kirwan type map
\begin{equation}\label{surj}\kappa : H_{K_0}^*(M, \Q) \to H_{K_0}^*(\mu^{-1}(a), \Q)\end{equation}
is surjective?}  

We will prove that the answer to this question is affirmative under certain restrictions.

\begin{proposition}\label{surjective} {\rm (Surjectivity criterium)} Assume that
all multiplicities $m_{\alpha}$ are strictly greater than 1 and  for any $b\in \a$, the set
$$Z_{\k}(b):=\{x\in \k : [x,b]=0\}=\k_0+\sum_{\alpha(b)=0}\k_{\alpha}$$
is the fixed point set of a certain torus $T_b\subset K_0$.  
Then the map $\kappa$ described by equation (\ref{surj}) is surjective, for any $a\in \a$.
\end{proposition}
\begin{proof} We consider the function $f: M\to \R$, $f(x)=\|\mu(x)-a\|^2$. By Theorem \ref{main} (i),
this is a minimally degenerate function. Moreover, $f$ is $K_0$-invariant. 
Let $C$ be a critical set of $f$ (by equation (\ref{refined}), $C$ can be $\mu^{-1}(a)$ or $C_{b,w}$). 
Denote $M^{\pm}=f^{-1}((-\infty, f(C)\pm\epsilon))$, for $\epsilon>0$ sufficiently small. 
By [K, chapter 10] (see also [BTW, section 9]) 
we have the commutative diagram
\begin{equation}\label{exactseq}
\vcenter{\xymatrix{
\cdots
\ar[r] &
H_{K_0}^*(M_+,M_-)
\ar[r]\ar[d]^\simeq &
H_{K_0}^*(M_+)
\ar[r] \ar[d] &
H_{K_0}^*(M_-)\ar[r] &  \cdots
\\
&
H_{K_0}^{*-{\rm index}(C) } (C)
\ar[r]^{~~ \ \  \ \ \ \ \cup e_C}  &
H_{K_0}^*(C) &
\\
}}
\end{equation}  
where $e_C\in H_{K_0}^*(C)$ denotes the equivariant Euler class of the normal bundle
$\nu(\Sigma_C)|_C$. 

We will prove that $e_C$ is not a divisor of zero. If $C=\mu^{-1}(a)$, then $e_C=1$ and the claim is obvious.
Let us consider  the case when $C=C_{b,w}=\mu^{-1}(a+b)\cap S_{wq,b}$ (see equation (\ref{refined})).
 According to a  criterium of Atiyah and Bott (see  [AB, Proposition 13.4]), it is
sufficient to prove that there exists a torus $T\subset K_0$ with the property that
the  only  points in $\nu(\Sigma_C)|_C$ which are fixed by $T$ are those from $C$. 
But $\nu(\Sigma_C)|_C$ is contained in $\nu(S_{wq,b})|_C$ (because $\Sigma_C$ contains
 $Y_{b,w}$ on a neighbourhood of $C$, and $\dim \Sigma_C=\dim Y_{b,w}$, see Theorem \ref{stratification}). We will show that the fixed points of 
the torus $T_b$ (see the statement of the proposition) on $\nu(S_{wq,b})$ are exactly those from
$S_{wq,b}$. First we note that the fixed point set of $T_b$ on $M$ is $$M\cap Z_{\p}(b)=Crit(h_b).$$
This is because if $\alpha$ is a positive root, then  the space $\p_{\alpha}$ is fixed by $T_b$ exactly
when $\k_{\alpha}$ is fixed by $T_b$, as $\k_{\alpha}=[x,\p_{\alpha}]$, for $x\in \a$ regular. 
But as we already mentioned in section  \ref{four}, $Crit(h_b)$ is the
disjoint union of all slices $S_{wq,b}$, where $w\in W$, and the conclusion follows.
\end{proof}

We will give an example where the criterium applies.

{\bf Example 1.} We consider the symmetric space $SU(2n)/Sp(n)$. For  details concerning this,
the reader can consult for instance [H, chapter X, section 2]. We will confine ourselves to saying
that the symplectic group 
$Sp(n)$ consists of all $n\times n$ nonsingular matrices with coefficients in $\H=\{a+ib+jc+kd: a,b,c,d\in \R\}$
which preserve the canonical symplectic product on $\H^n$. There is a canonical embedding of this group
in $SU(2n)$, as the fixed point set of a certain involutive automorphism of the group $SU(2n)$.    
In this case, we can choose $\a$ to be the space  of all matrices of  type
$$a:=\left(%
\begin{array}{ccccc}
   iD&  0 \\
  0 & iD
\end{array}%
\right),
$$ 
where $i=\sqrt{-1}$ and  $$D=\left(%
\begin{array}{ccccc}
  a_1 & 0 &\ldots & 0 &0\\
  0 & a_2 & 0 & \ldots & 0  \\
  \vdots & \vdots & \vdots &\vdots & \vdots \\
  0 & 0 & 0 & \ldots & a_n
\end{array}%
\right),
$$ 
with $a_1, \ldots, a_n \in \R$, $a_1 +a_2 +\ldots +a_n=0$. The positive roots are $\alpha_{rs}:\a\to \R$,
$$\alpha_{rs}(a)=a_r-a_s,$$
where $r< s$. This is an irreducible root system of type $A_{n-1}$.
  We have $$\k_{\alpha_{rs}}=\{\gamma E_{rs}-\bar{\gamma}E_{sr} : \gamma \in \H\},$$
where $E_{rs}$ denotes the $n\times n$ matrix whose entries are all 0, except the one on line $r$ and
column $s$, which is 1. In particular, all multiplicities are equal to $4$. 
The space $K_0$ is just $Sp(1)^n$. We can describe explicitely the action of an arbitrary element
$h=(h_1,\ldots, h_n)\in Sp(1)^n$ on $\k_{\alpha_{rs}}$, as follows:
\begin{equation}\label{E}h.(\gamma E_{rs}-\bar{\gamma}E_{sr}) = h_r\gamma \bar{h}_s E_{rs}-h_s\bar{\gamma}\bar{h}_r E_{sr}.
\end{equation}

We show that the fixed point set of the group $$Sp(1)_{rs}:=\{(h_1,\ldots,h_{r-1}, 1,h_{r+1},\ldots,
h_{s-1},1, h_{s+1},\ldots , h_n)\}\subset K_0$$ on $\k$  is $\k_0+\k_{\alpha_{rs}}$. This follows from
equation (\ref{E}) and the following technical claim.

{\it Claim.} {\it (a) If $\gamma \in \H$ such that $h_1\gamma \bar{h}_2=\gamma$ for any $h_1, h_2\in \H$ of length 1,
then $\gamma =0$.

(b) The same is true if $h_1, h_2 $ are of the form $a+ib$, where $a,b\in \R$, $a^2+b^2=1$.}

The proof of the claim is straightforward. In fact, point (b) shows that  $\k_0+\k_{\alpha_{rs}}$ is
the fixed point set of the torus
$$T_{rs}:=\{(h_1,\ldots,h_{r-1}, 1,h_{r+1},\ldots,
h_{s-1},1, h_{s+1},\ldots , h_n): h_p\in \bC, |h_p|=1\}\subset K_0$$
In a similar way, we can prove that the hypothesis of Proposition \ref{surjective} is satisfied for
any $b\in \a$.

Now we give another example, where the criterium does not apply.

{\bf Example 2.} Consider the symmetric space $\bC P^2 = SU(3)/U(2)$, where $U(2)$ is embedded in $SU(3)$ via
$$A\mapsto \left(%
\begin{array}{ccccc}
  \frac{1}{\det(A)} & 0 \\
  0 & A 
\end{array}%
\right),
$$ 
$A\in U(2)$. It turns out that $\p$ consists of all matrices of the type
$$\left(%
\begin{array}{ccccc}
  0 & -\bar{z}_1 & -\bar{z}_2 \\
  z_1 & 0 & 0  \\
  z_2 & 0 & 0 
\end{array}%
\right),  $$
where $z_1, z_2 \in \bC$. We choose $\a$ the space of all
matrices $$a:=\left(%
\begin{array}{ccccc}
  0 & -{x} & 0 \\
  x & 0 & 0  \\
  0 & 0 & 0 
\end{array}%
\right),$$  with $ x \in \R$. The  positive roots are $\alpha$ and $2\alpha$, where
$\alpha(a)=-x$ and the corresponding root spaces in $\p$ are
$$\p_{\alpha}= \{\left(%
\begin{array}{ccccc}
  0 & 0 & -\bar{z} \\
  0 & 0 & 0  \\
  z & 0 & 0 
\end{array}%
\right): z\in \bC\}, \quad \p_{2\alpha}=\{\left(%
\begin{array}{ccccc}
  0 & iy & 0 \\
  iy & 0 & 0  \\
  0& 0 & 0 
\end{array}%
\right): y\in \R\}.$$
We deduce that $m_{\alpha}=2+1=3$. An easy calculation shows that $K_0:=Z_K(\a)$ is the subgroup of $U(2)$ which
consists of 
$$\left(%
\begin{array}{ccccc}
  z & 0 & 0 \\
  0 & z & 0  \\
  0& 0 & \frac{1}{z^2} 
\end{array}%
\right)$$
for $z\in \bC\setminus \{0\}$. One can see that $K_0$ acts trivially  not only on $\a$, but also  
on $\p_{2\alpha}$. We deduce that it acts trivially on $\k_{2\alpha}$ as well. This implies
that there is no subgroup of $K_0$ whose fixed point set in $\k$ is just $\k_0$.
 So the hypothesis of Proposition \ref{surjective} is not satisfied.

\end{document}